\documentclass[12pt]{article}
\usepackage{epsfig}
\usepackage{amsmath}
\usepackage{amssymb}
\usepackage{amsfonts}
\usepackage{amsthm}

\newcommand{\R}{\mathbb{R}}
\newcommand{\eps}{\varepsilon}
\newcommand{\loc}{\textrm{loc}}
\newcommand{\br}{\begin{eqnarray}}
\newcommand{\er}{\end{eqnarray}}
\newcommand{\no}{\nonumber}

\newtheorem{lemma}{Lemma}[section]
\newtheorem{theorem}[lemma]{Theorem}
\newtheorem{corollary}[lemma]{Corollary}
\newtheorem{proposition}[lemma]{Proposition}

\newtheorem{definition}[lemma]{Definition}

\newcommand{\nit}{\noindent}
\newcommand{\be}{\begin{equation}}
\newcommand{\ee}{\end{equation}}

\newcommand{\lam}{\mbox{$\lambda$}}

\newcommand{\abs}[1]{\lvert #1 \rvert}
\newcommand{\norm}[1]{\lVert #1 \rVert}

\begin{document}

\title{Existence of KPP fronts in spatially-temporally periodic advection and variational principle for propagation speeds}

\author{James Nolen\thanks{Department of Mathematics, University of Texas at Austin,
Austin, TX 78712 (jnolen@math.utexas.edu).}
\and Matthew Rudd\thanks{Department of Mathematics, University of Texas at Austin,
Austin, TX 78712 (rudd@math.utexas.edu).}
\and Jack Xin\thanks{Department of Mathematics and ICES (Institute of Computational
Engineering and Sciences), University of Texas at Austin, Austin, TX 78712
(jxin@math.utexas.edu).}}

\date{}
\setcounter{page}{1}
\maketitle
\begin{abstract}
We prove the existence of Kolmogorov-Petrovsky-Piskunov (KPP) type
traveling fronts
in space-time periodic and mean zero incompressible advection, and
establish a variational (minimization) formula for the minimal speeds.
We approach the existence
by considering limit of a sequence of front solutions to
a regularized traveling
front equation where the nonlinearity is
combustion type with ignition cut-off.
The limiting front equation is
degenerate parabolic and does not permit strong solutions, however, the
necessary
compactness follows from monotonicity of fronts and degenerate regularity.
We apply a dynamic argument to justify that the constructed KPP
traveling fronts
propagate at minimal speeds, and derive the speed variational formula.
The dynamic method avoids the degeneracy in traveling front equations, and
utilizes the parabolic maximum principle of the governing
reaction-diffusion-advection equation.
The dynamic method does not
rely on existence of traveling fronts.

\end{abstract}

\newpage

\begin{section}{Introduction} \label{intro}
\setcounter{equation}{0}

%\makeatletter
%\renewcommand{\theenumi}{\roman{enumi}}
%\renewcommand{\labelenumi}{(\theenumi)}
%\makeatother

Reaction-diffusion front propagation in heterogeneous fluid flows
has been an active research area for decades
(\cite{CW}, \cite{Kh}, \cite{KA}, \cite{Ro},
\cite{Vlad}, \cite{Xin1}, \cite{Yak} and references therein).
A fascinating phenomenon is that the large time (large scale) front speed
is a highly nontrivial quantity that depends on the multiple scale
structures of the fluid flows.
Speed characterizations have been studied mathematically
for various flow patterns by analysis of
the reaction-diffusion-advection (RDA) equations (see
\cite{BH1,Const1,Ham99,HPS,KR,Kh,MS1,MS2,NX1,PX,Vlad,Xin1,Xin2} and references
therein). In particular, variational characterizations have led to speed
asymptotics in both deterministic and random flows \cite{PX,HPS,Xin2,nolen:ekt04,NX3,NX4}.

In this paper, we consider front solutions to a RDA equation of the form:
\begin{equation} \label{orig-eqn}
u_{t}  =  \Delta u + b(x,t) \cdot \nabla u + f(u) \, ,
\end{equation}
where $u = u(x,t)$, $x \in \R^{N}$, and $t \in \R$.
The $N$ components of the vector field $ b(x,t) := ( b^{1}(x,t), b^{2}(x,t), \ldots, b^{N}(x,t) ) $
are smooth and spatially divergence--free, are periodic of period $1$ in both $x$ and $t$, and
have mean zero over the period cell $Q \times (0,1)$, where $Q$ is the unit cube in $\R^{N}$:
\begin{equation} \label{flow-field}
\nabla_{x} \cdot b = 0 \quad \textrm{and} \quad
\int_{Q \times (0,1)}{ b^{i}(x,t) \, dx \, dt } = 0, \quad \textrm{for} \quad i = 1, \ldots, N \, .
\end{equation}
The function $f : [0,1] \rightarrow \R$ is KPP, namely a $C^2$ function of $u$
such that:
\begin{equation} \label{kpp-f}
f(0) = f(1) = 0, \; f'(0) > 0, \; f'(1) < 0; \;
f(u) > 0, \forall  u \in (0,1); \;f(u) \leq uf'(0).
\end{equation}
An example is $f(u)=u(1-u)$.

A more general class of nonnegative reaction nonlinearities arising
in applications may
allow degeneracy at zero, e.g. $f(u) = u^m(1-u)$, for $m \geq 2$,
and $f(u) = e^{-E/u}(1-u), E>0$. Let us define the so called
positive nonlinearity to include KPP and generalizations. A {\it
positive nonlinearity} $f(u)$ satisfies: $f(0)=f(1)=0$; $f(u) > 0$,
$u \in (0,1)$.
In the course of this study we will also
consider the \textit{ignition type} combustion nonlinearity,
for which $f(1) = 0$, $f'(1) < 0$, and there exists a constant $\theta \in (0,1)$ such that
\begin{equation} \label{combustion-f}
f(u) = 0 \quad \textrm{for} \quad u \in [ 0, \theta ] \quad \textrm{and} \quad f(u) > 0 \quad \textrm{for} \quad u \in (\theta,1) \, .
\end{equation}
The constant $\theta$ is known as the \textit{ignition temperature}.

Given a unit vector $k = (k^{1},\ldots,k^{N}) \in \R^{N}$,
we look for a planar
traveling front which propagates in the
direction $k$ with constant speed and solves (\ref{orig-eqn}).
Due to the form of the advection field, we seek a solution of the form
\begin{equation} \label{ansatz}
u(x,t) = U(k \cdot x - ct, x, t) = U(s,y,\tau)
\end{equation}
for $s := k \cdot x - ct$, $y := x$, and $\tau := t $. The constant $c$ is the speed of the front.
We require $U(s,y,\tau)$ to be 1--periodic in both $y$ and $\tau$ and
to satisfy the following conditions at $s = \pm \infty$,
uniformly in $y$ and $\tau$:
\begin{equation}
\quad \lim_{s \to -\infty} U(s,y,\tau) = 1
\quad \text{and} \quad \lim_{s \to \infty} U(s,y,\tau) = 0 \, .
\end{equation}
Substituting $U$ into equation (\ref{orig-eqn}) and
imposing these auxiliary conditions yields the problem
\begin{equation} \label{U}
\left\{ \begin{array}{c}
U_{\tau}  =  \left( k \partial_{s} + \nabla_{y} \right)^{2} U +
 b \cdot \left( k \partial_{s} + \nabla_{y} \right) U + c \, U_{s} + f( U ) \, , \\
\\
\displaystyle{ U(s,y,\tau) \quad 1\textrm{--periodic in} \; y \; \textrm{and}
\;  \tau, } \\
\\
\displaystyle{ \lim_{s \to \infty} U(s,y,\tau) = 0 \quad \text{and}
\quad \lim_{s \to -\infty} U(s,y,\tau) = 1 } \, .
\end{array} \right.
\end{equation}
Equation (\ref{U}) contains only directional derivatives that do not span
$\nabla_{s,y,\tau}\, U$, hence there is a lack of compactness. One can view (\ref{U})
as a degenerate parabolic equation, parabolic along $(s,y,\tau) = (-c,0,1)$,
and degenerate elliptic in $(s,y,\tau)=(k^{i},e^{i},0)$, $e^{i}$ the standard i-th
unit vector in $R^N$, $i=1,2,\cdots,N$.
To proceed with existence,
we shall regularize (\ref{U}) then pass to the limit to prove a certain weak solution
of (\ref{U}) defined below.

\begin{definition} \label{front-defn}
A \textit{front traveling at speed $c \neq 0$ with reaction $f$}
is a locally integrable function $U(s,y,\tau)
\in L^{1}_{\loc}(\R \times Q \times (0,1))$ whose directional derivatives
\[
U_{\tau} - c U_{s}, \quad k^{i}\, U_{s} + U_{y^{i}},
\ i = 1, \ldots, N, \quad \textrm{and}
\quad \left( k \partial_{s} + \nabla_{y} \right)^{2} U
\]
are continuous and satisfy the traveling front equation (\ref{U}) in the form:
\begin{equation} \label{front-eqn}
U_{\tau} - c U_{s}  =  \left( k \partial_{s} + \nabla_{y} \right)^{2} U +
b \cdot \left( k \partial_{s} + \nabla_{y} \right) U + f( U ) \, ,
\end{equation}
and the associated boundary conditions.
\end{definition}

Our first main result is that such a front solution exists for
positive nonlinearity.
The proof is based on elliptically
regularizing (\ref{front-eqn}) and approximating
positive nonlinearity by a sequence of combustion
nonlinearities with ignition cut-off.
The regularized problem can be solved by existing methods in the literature,
and one key property is that $U$ is monotone in $s$.
The monotonicity, the analysis of
asymptotics near $s$ infinities, and the available derivative estimates of
(\ref{front-eqn}) imply enough compactness for passage of limit and
constructing a desired solution.

A few remarks are in order on the traveling fronts.
In the original coordinates $(x,t) \in \R^{N} \times \R$,
we see that a traveling front $U$ corresponds to the
function $u(x,t)=U(k \cdot x - ct, x, t)$. This is a classical solution of
(\ref{orig-eqn}), because derivatives in the original $(x,t)$ coordinates
translate to directional derivatives in traveling front variables $(s, y, \tau)$.
Note also that this definition of the traveling front is
different from the definition of pulsed traveling fronts in
time-independent advection fields \cite{BH1}.
In that context, a pulsed traveling front
in direction $k \in S^{N-1}$ is defined by the relation
\[
u(x + e,t) = u(x,t - \frac{k \cdot e}{c}), \;\;\forall (x,t) \in R^{N+1},
\]
for any vector $e \in R^N$ such that $b(x + me) \equiv b(x)$ for some integer $m$.
If this relation were imposed on the solutions considered here,
we would deduce that
\br
u(x + e,t) &=& U(k\cdot x + k\cdot e - ct, x + e, t) \no \\
& = & U(k\cdot x - c (t - \frac{k\cdot e}{c}), x, t) \no
\er
would be equal to
\[
u(x, t- \frac{k\cdot e}{c}) =
U(k \cdot x - c(t - \frac{k\cdot e}{c}), x, t- \frac{k\cdot e}{c}).
\]
Hence, if $b$ is $T$-periodic in time variable $t$, we would conclude that
\[
c = \frac{k\cdot e}{mT},\;\;\;\text{for some integer}\; m\neq 0
\]
This implies that if $b$ has a nontrivial time-dependence,
the speed $c$ would be determined only by the period of $b$ in space and time,
regardless of other information on $b$, which is nearly impossible.
The price of seeking a more general form of traveling fronts (\ref{ansatz})
is the difficulty caused by degeneracy in front equation (\ref{front-eqn}).

Degeneracy makes it hard to prove the minimality of front speed,
if one relies entirely on the front equation (\ref{front-eqn}).
We adopt the alternative
approach of studying the time dependent
front-like solutions of (\ref{orig-eqn})
where {\it there is dynamics but no degeneracy}. The KPP nonlinearity
is critically used in constructing comparison functions and passing to the
large time limits of solutions.

Our second main result is that the KPP
traveling front speed is minimal (denoted by $c^*$) and is given by
the variational principle:
\[
c^* = \inf_{\lam > 0}\, \mu (\lam) / \lam,
\]
where $\mu (\lam)$ is the
principal eigenvalue of the periodic--parabolic operator \cite{hess:ppb91}
\[
L^{\lambda} \Phi :=
\Delta_{x} \Phi + \left( b - 2 \lambda k \right) \cdot \nabla_{x}\Phi +
\left( \lambda^{2} - \lambda b \cdot k + f'(0) \right) \Phi - \Phi_{t},
\]
defined on spatially--temporally periodic functions $\Phi(x,t)$.
The variational principle reveals that the flow structures in $b$ are
upscaled into the front speed $c^*$ via the principal eigenvalue $\mu$, 
see \cite{GFr79} for a similar formula.

The rest of the paper is organized as follows. In section 2, we study
regularized traveling fronts under combustion nonlinearity
with cut-off, and obtain necessary estimates to remove the regularization, and
justify boundary conditions at $s$ infinities. Front existence for
positive nonlinearities follows. In section 3, we study
propagation of KPP fronts, minimality of front speeds and their variational
characterization. In section 4, we present properties of the
associated eigenvalue problem. Concluding remarks are in section 5.

\end{section}

\begin{section}{Existence of Fronts: Positive Nonlinearity} \label{section-existence-comb}
\setcounter{equation}{0}

In this section, we study the problem of existence of planar traveling front solutions of (\ref{orig-eqn})
for general positive nonlinearities. We prove the following:

\begin{theorem} \label{existence-kpp}
Let $f$ be a positive nonlinearity: $f(u)> 0$, $u\in (0,1)$; $f(0)=f(1)=0$.
Then there exists a traveling front $U^*$
with speed $c^{*} > 0$ that solves (\ref{front-eqn}) with the
boundary conditions in
(\ref{U}).
\end{theorem}

\noindent
Our proof of Theorem \ref{existence-kpp} proceeds by first solving a regularized version of (\ref{U}) where
we approximate $f$ by a combustion-type nonlinearity $f^\theta$
such that $f^\theta \to f$ as $\theta \to 0$. In Section \ref{regularized-combustion}, we prove existence of solutions to the regularized problem for the combustion nonlinearity $f^\theta$. Then, in Section \ref{regularized-kpp} we show that as $\theta \to 0$, $f^\theta \to f$ and the solutions converge to solutions of the regularized problem for the general positive nonlinearity. Finally, in Section \ref{degenerate-kpp} we complete the proof of Theorem \ref{existence-kpp} by removing the regularization parameter to recover
a desired solution.

\begin{subsection}{Regularized Equation: Combustion Nonlinearity} \label{regularized-combustion}
In this section, we assume $f = f^\theta$ is a combustion-type nonlinearity, $f(u) = 0$ if $u \in [0, \theta]$. Let $\eps > 0$ be given, and consider the family of regularized problems
\begin{equation} \label{param-eps}
\left\{ \begin{array}{c}
\displaystyle{ U^{\sigma}_{\tau} =  \eps U^{\sigma}_{ss} + \left( k \partial_{s} + \nabla_{y} \right)^{2} U^{\sigma}
+ \sigma b \cdot \left( k \partial_{s} + \nabla_{y} \right) U^{\sigma} + c_{\sigma} \, U^{\sigma}_{s} + f( U^{\sigma} ) \, , } \\
\\
\displaystyle{ U^{\sigma}(s,y,\tau) \quad 1\textrm{--periodic in} \; x \; \textrm{and} \;  \tau, } \\
\\
\displaystyle{ \lim_{s \to \infty} U^{\sigma}(s,y,\tau) = 0 \, , \quad \textrm{and} \quad \lim_{s \to -\infty} U^{\sigma}(s,y,\tau) = 1 } \, ,
\\
\end{array} \right.
\end{equation}
where the parameter $\sigma \in [0,1]$. In this section, we will usually suppress the dependence of $U^\sigma = U^{\sigma, \epsilon, \theta}$  and $c_\sigma = c_{\sigma, \epsilon, \theta}$ on $\epsilon$ and on $\theta$ since they remain fixed throughout the section.

Before discussing the existence of solutions of (\ref{param-eps}), let us point out some of the
properties that solutions of (\ref{param-eps}) must possess.  Since the equation in (\ref{param-eps}) is nondegenerate, one may directly use the
classical parabolic maximum principle to prove the following (cf. \cite{xin:epf92},\cite{nolen:ekt04}):
\begin{itemize}

\item
$ 0 < U^{\sigma} < 1 \quad \textrm{on} \quad \R \times Q \times (0,1)$.

\item
$ U^{\sigma}(s,y,\tau) > U^{\sigma}(s',y,\tau)$ for $s < s'$ and $(y,\tau) \in Q \times (0,1)$ (monotonicity).

\item
The speed $c_{\sigma}> 0$ for which (\ref{param-eps}) has a solution $(U^{\sigma}, c_{\sigma})$ must be unique.

\item
The solution $U^{\sigma}$ must be unique modulo translation in $s$.

\item
The solution $U^\sigma$ converges exponentially fast to the limits as $s \to \pm \infty$.
\end{itemize}

Our goal is to find a classical solution for this regularized problem when $\sigma = 1$. It is known (\cite{xin:epf92},\cite{nolen:ekt04}) that there exists a solution of (\ref{param-eps}) when $\sigma = 0$ (i.e., in the absence of advection), so we apply a continuation argument to perturb this solution and solve (\ref{param-eps}) when $\sigma = 1$. We now outline the application of this method to (\ref{param-eps}); the argument generalizes the arguments of \cite{xin:epf92} and \cite{nolen:ekt04}.

Assuming that (\ref{param-eps}) has a classical solution $(U^{0},c_{0})$ when
$\sigma = \sigma_{0}$, we seek a solution $ ( U^{\sigma}, c_{\sigma} ) $ of (\ref{param-eps}) for small perturbations of $\sigma > \sigma_{0}$.  Substituting the expressions
\[
c_{\sigma} = c_{0} + \sigma \bar{c}_{\sigma} \quad \textrm{and} \quad
U^{\sigma} = U^{0} + \sigma V^{\sigma}
\]
into equation (\ref{param-eps}) and using the fact that
\[
U^{0}_{\tau}  =
\eps U^{0}_{ss} + \left( k \partial_{s} + \nabla_{y} \right)^2 U^{0} + \sigma_{0} b \cdot \left( k \partial_{s} + \nabla_{y} \right) U^{0}
+ c_{0} U^{0}_{s} + f( U^{0} ) \, ,
\]
we obtain the following nonlinear fixed-point problem for the perturbations $\tilde{c}_{\sigma}$ and $V^{\sigma}$:
%% We first obtain
%% \[
%% \begin{array}{rcl}
%% U^{0}_{\tau} + \sigma V^{\sigma}_{t} & = &
%% U^{0}_{ss} + 2 \, k^{i} U^{0}_{s, x^{i}} + \Delta_{x} U^{0} + ( \sigma b \cdot k + c_{0} + \sigma c_{\sigma} ) \, U^{0}_{s} + \sigma b \cdot \nabla_{y} U^{0} \\
%% & & \\
%% & &
%% + \sigma \left( V^{\sigma}_{ss} + 2 \, k^{i} V^{\sigma}_{s, x^{i}} + \Delta_{x} V^{\sigma} + ( \sigma b \cdot k + c_{0} + \sigma c_{\sigma} ) \, V^{\sigma}_{s}
%% + \sigma b \cdot \nabla_{y} V^{\sigma} \right) \\
%% & & \\
%% & &
%% + f( U^{\sigma} ) \, ,
%% \end{array}
%% \]
%% dividing through by $\sigma$, and subtracting $f'(U^{0}) V^{\sigma}$ from both sides, we have
\begin{equation} \label{v-sigma}
V^{\sigma}_{\tau} - \eps V^{\sigma}_{ss} - \left( k \partial_{s} + \nabla_{y} \right)^{2} V^{\sigma} - c_{0} V^{\sigma}_{s}
- f'(U^{0}) V^{\sigma} = F( V^{\sigma}, \sigma, \bar{c}_{\sigma} ) \, ,
%% \left( b \cdot k + c_{\sigma} \right) U^{0}_{s} + b \cdot \nabla_{y} U^{0}
%% + ( \sigma b \cdot k + \sigma c_{\sigma} ) \, V^{\sigma}_{s} \\
%% + \sigma b \cdot \nabla_{y} V^{\sigma} + \left[ \frac{ f( U^{\sigma} ) - f( U^{0} ) }{ \sigma } - f'(U^{0})V^{\sigma} \right] \, .
\end{equation}
where the terms not appearing explicitly on the left have been absorbed into the nonlinearity $F$.

To establish the existence of a solution $V^{\sigma}$ of (\ref{v-sigma}), we employ the arguments used in \cite{xin:epf92}.
Consider the weighted norm
\begin{equation} \label{weight-norm}
\| u \|_{\rho}^{2} := \int_{\R \times Q \times (0,1)}{ u^{2} \, \rho \, ds \, dy \, d \tau } ,
\end{equation}
where the weight $\rho : \R \rightarrow (0,\infty)$ is defined by
\[
\rho(s) := 1 + e^{-2 \delta s }
\]
for an appropriate constant $0 < \delta < 1$.  The weighted spaces
$L^{2}_{\rho} := L^{2}_{\rho}(\R \times Q \times (0,1))$ and $H^{1}_{\rho} : = H^{1}_{\rho}(\R \times Q \times (0,1) )$
are defined in the natural way.  Define the simple function $A : \R \rightarrow \R$ by
\[
A(s) := \left\{ \begin{array}{rcl}
f'(1), & \textrm{if} & s > 0, \\
& & \\
0, & \textrm{if} & s \leq 0,
\end{array} \right.
\]
and let $L$ and $L_{1}$ be operators with the common domain
\[
D(L) = D( L_{1} ) := \left\{ \, v \in H^{1}_{\rho}(\R \times Q \times (0,1)) \, : \,
\left( k \partial_{s} + \nabla_{y} \right)^{2} v \in L^{2}_{\rho} \, \right\}
\]
and respective definitions
\begin{equation} \label{weight-op}
L v := v_{\tau} - \eps v_{ss} - \left( k \partial_{s} + \nabla_{y} \right)^{2} v - c_{0} v_{s} - f'(U^{0}) v \quad \textrm{and}
\end{equation}
\begin{equation} \label{weight-op1}
L_{1} v := v_{\tau} - \eps v_{ss} - \left( k \partial_{s} + \nabla_{y} \right)^{2} v - c_{0} v_{s} - A(s) v \, ,
\end{equation}
for $v \in D(L)$. Observe that $L$ is the linear operator defined by the left-hand side of equation (\ref{v-sigma}) and that we want to solve $L V^\sigma = F(V^\sigma,\sigma, \bar c_\sigma)$.

As in \cite{nolen:ekt04} and \cite{xin:epf92}, we now consider the invertibility of the operator $L$. Because of the regularization with $\epsilon>0$, we can use the maximum principle and standard parabolic estimates to show the following:
$L_{1}$ is invertible on $L^{2}_{\rho}$, $L_{1}$ and $L$ differ by a relatively compact operator, and
$L$ is a Fredholm operator of index zero, whose adjoint $L^{*}$ has a $1$--dimensional null space and an isolated eigenvalue
of finite multiplicity at $0$.
Moreover, the eigenfunction $v^{*}$ of $L^{*}$ corresponding to $0$ is strictly positive,
and the following Fredholm alternative holds:
\begin{proposition} \label{fat}
Let $f \in L^{2}_{\rho}$.  The equation $L v = f$ has a solution if and only if
\begin{equation} \label{fat-eqn}
\int_{\R \times Q \times (0,1)}{ f v^{*} \, \rho \, ds \, dy \, d\tau } = 0 \, .
\end{equation}
When the solvability condition (\ref{fat-eqn}) holds, the solutions of $ L v = f $ form a one--dimensional space.
\end{proposition}

\noindent
This result is a straightforward generalization of the corresponding results in \cite{nolen:ekt04} and \cite{xin:epf92}, so we do not include the details.  Now, to solve equation (\ref{v-sigma}), we can use Proposition \ref{fat} to set up an iteration scheme, as in Section 2.3 of \cite{xin:est91}, to which the contraction mapping principle applies.
With each iteration, the constants $\{\bar c^n\}_{n=1}^\infty$ that converge to $\bar c_\sigma$ are determined by the solvability condition (\ref{fat-eqn}). As a result, we find that (\ref{param-eps}) is solvable for small perturbations of $\sigma$:
\begin{theorem} \label{continue}
If problem (\ref{param-eps}) has a classical solution $(U^{0}, c_{0})$ when $\sigma = \sigma_{0}$,
then there exists a positive constant $\sigma^{*} = \sigma^{*}(U^{0}, c_{0})$ such that (\ref{param-eps}) has a unique classical
solution $ ( U^{\sigma}, c_{\sigma} ) $ whenever $ \sigma_{0} < \sigma < \sigma^{*} $.
\end{theorem}
%% As noted above, uniqueness follows from the normalization condition.
\noindent For details of this iteration procedure, see \cite{nolen:ekt04} and \cite{xin:est91}.

Next, having shown that solutions $(U^{\sigma},c_{\sigma})$ of (\ref{param-eps}) exist for $\sigma < \sigma^{*}$, we consider
the behavior of solutions as $\sigma \rightarrow \sigma^{*}$. We first prove bounds on the sequence of speeds $\{ c_\sigma \}$ so that we can extract a convergent subsequence. Using periodicity and the conditions at infinity, integrating
equation (\ref{param-eps}) over $\R \times Q \times (0,1)$ yields
\begin{equation} \label{speed-1}
 c_{\sigma} = \int{ f( U^{\sigma} ) \, ds \, dy \, d\tau } > 0 \, ,
\end{equation}
from which we see that all of the speeds $ c_{\sigma} $ are positive.
Furthermore, as in \cite{heinze:hff}, multiplying equation (\ref{param-eps}) by $U^{\sigma}$ and
integrating over $\R \times Q \times (0,1)$ yields the $L^{2}$ estimate
\begin{equation} \label{L2-dsdx}
\int{ \left| k U^{\sigma}_{s} + \nabla_{y} U^{\sigma} \right|^{2} \, ds \, dy \, d\tau } =
-\frac{ c_{\sigma} }{ 2 } + \int{ U^{\sigma} f( U^{\sigma} ) \, ds \, dy \, d\tau } \leq
 \frac{c_{\sigma}}{2} \, ,
\end{equation}
where we used the identity (\ref{speed-1}) and the fact that $0 < U^{\sigma} < 1$ to bound the middle term.
In addition, we have the following essential estimates on the speeds.
\begin{proposition} \label{wave-speed}
There exist constants $\underline{c}$ and $\overline{c}$, independent of both $\sigma$ and $\eps$, such that
\begin{equation} \label{wave-speed-ineq}
0 < \underline{c} \, \leq \, c_{\sigma, \epsilon} \, \leq \, \overline{c} \quad \textrm{for all} \quad \sigma \in [0,\sigma^*) ,\, \epsilon \in (0,1]\, .
\end{equation}
\end{proposition}

\begin{proof}

If $ ( \, V, \, \overline{c}_\epsilon \, ) $ solves the one-dimensional problem
\begin{equation} \label{1d-fast}
(1 + \epsilon) V_{ss} + \left( \alpha + \overline{c}_\epsilon \right) V_{s} + f(V) = 0 \, ,
\end{equation}
with the appropriate limits at $s = \pm \infty$ and the constant coefficient $\alpha$ defined by
\[
\alpha := \inf_{Q \times [0,1]}{ \sigma^* b(y,\tau) \cdot k } \leq 0 \, ,
\]
then the proof of Lemma 3.1 in \cite{xin:epf92} shows that
\[
c_{\sigma,\epsilon}  \leq  \overline{c}_\epsilon , \quad \textrm{for all} \quad \sigma \in [0,\sigma^*) . 
\]
Note that since $\alpha$ a constant, the constant 
\be
\tilde c_\epsilon = \frac{\alpha + \overline{c}_\epsilon}{(1 + \epsilon)} \label{tildeeps}
\ee
is the unique speed associated with the one-dimensional traveling wave problem 
\begin{equation} 
V_{ss} + \tilde c_\epsilon V_{s} + \frac{1}{(1 + \epsilon)} f(V) = 0 \, . \label{tildepseq}
\end{equation}
(See \cite{berestycki:tfc92} for proof of existence and uniqueness of $\tilde c_\epsilon$.) Because $\epsilon \in (0,1]$, a comparison argument as in \cite{xin:ijm91} shows that $\tilde c_\epsilon$ is bounded above by $\tilde c_0 > 0$, the unique speed determined by (\ref{tildepseq}) with $\epsilon = 0$. Therefore, from (\ref{tildeeps}) we have
\[
\overline{c}_\epsilon \leq (1 + \epsilon) \tilde c_0 - \alpha \leq 2 \tilde c_0 + \sigma^* \norm{b}_\infty = \overline{c}, \;\; \forall\, \epsilon \in (0,1].
\]
This proves the upper bound in (\ref{wave-speed-ineq}).
As in \cite{heinze:hff}, define the constant $\underline{c}$ by
\[
\underline{c} := \int_{0}^{1}{ \sqrt{ 2 f(\xi) } \, d\xi } \, .
\]
We have
\[
\begin{array}{rcl}
\displaystyle{ \left( \int{ \sqrt{ 2 f(\xi) } \, d\xi } \right)^{2} } & = &
\displaystyle{ \left( \int{ \sqrt{ 2 f( U^{\sigma} ) } \, k \cdot \left( k U^{\sigma}_{s} + \nabla_{y} U^{\sigma} \right) \, ds \, dy \, d\tau } \right)^{2} } \\
& & \\
& \leq &
\displaystyle{ \left( \int{ f( U^{\sigma} ) \, ds \, dy \, d\tau } \right) \left(  2 \int{ \left| k U^{\sigma} + \nabla_{y} U^{\sigma} \right|^{2} \, ds \, dy \, d\tau }
\right) } \\
& & \\
& \leq &
( c_{\sigma} )^{2} \, ,
\end{array}
\]
from which the lower bound follows. Here we have used (\ref{L2-dsdx}) and the fact that $U^\sigma_s < 0$ for all $(s,y,\tau) \in R \times Q \times [0,1]$.

\end{proof}

Proposition \ref{wave-speed} proves that the wave speeds are bounded away from both $0$ and $+\infty$ and
thus have a positive limit point as $\sigma \rightarrow \sigma^{*}$.
Using this fact and standard parabolic regularity results, we obtain the following:
\begin{theorem} \label{existence-eps}
Let $f$ be a combustion nonlinearity.  For each $\eps > 0$, there exists a classical solution $(U^{\eps},c_{\eps})$ of
\begin{equation} \label{eqn-eps}
\left\{ \begin{array}{c}
\displaystyle{ U^{\eps}_{\tau} =  \eps U^{\eps}_{ss} + \left( k \partial_{s} + \nabla_{y} \right)^{2} U^{\eps}
+ b \cdot \left( k \partial_{s} + \nabla_{y} \right) U^{\eps} + c_{\eps} \, U^{\eps}_{s} + f( U^{\eps} ) \, , } \\
\\
\displaystyle{ U^{\eps}(s,y,\tau) \quad 1\textrm{--periodic in} \quad y \quad \textrm{and} \quad  \tau, } \\
\\
\displaystyle{ \lim_{s \to \infty} U^{\eps}(s,y,\tau) = 0 \, , \quad \textrm{and} \quad \lim_{s \to -\infty} U^{\eps}(s,y,\tau) = 1 } \, .
\\
\end{array} \right.
\end{equation}
The solution $U^{\eps}$ is unique modulo translation in $s$.  Moreover,
\[
0 < U^{\eps} < 1 \quad \textrm{on} \quad \R \times Q \times [0,1] \, ,
\]
\[
U^{\eps}(s,y,\tau) > U^{\eps}(s',y,\tau) \quad \textrm{for} \quad s < s', \ y \in Q, \ t \in [0,1] \, , \quad \textrm{and}
\]
\[
c_{\eps} > 0 \, .
\]
\end{theorem}

\begin{proof}

Fix an $\eps > 0$.  By interior parabolic regularity results (\cite{lieberman:sop96}) and the bounds on $f$ and $U^{\sigma}$,
\be
\| D^{2} U^{\sigma} \|_{p,E} + \| U^{\sigma}_{\tau} \|_{p,E} \leq
C \left( \| U^{\sigma} \|_{p,E} + \| f( U^{\sigma} ) \|_{p,E} \right)
\leq C \left( 1 + \| f \|_{\infty} \right) | E | \, , \label{wtwopest}
\ee
where $E \subset \subset \R \times Q \times (0,1)$ and the constant $C$ is independent of $\sigma$.
As the upper bound on the right is independent of $\sigma$, choosing $p$ sufficiently large ensures
that $U^{\sigma}$ is H\"{o}lder continuous.  Schauder estimates then apply, and the
Arzel\'{a}--Ascoli Theorem  guarantees that $U^{\sigma} \rightarrow U^{*}$
uniformly on any compact set as $\sigma \rightarrow \sigma^{*}$, where $U^{*}$ solves the
equation in (\ref{param-eps}) for $\sigma^{*}$ (here we pass to a convergent subsequence, if necessary).  Since problem (\ref{param-eps}) has a solution
at $\sigma^{*}$, solvability may actually be continued all the way to $\sigma = 1$, yielding
the solution $U^{\eps}$ of problem (\ref{eqn-eps}).
Because the convergence of the functions $U^\sigma$ is uniform, the boundary conditions at $s = \pm \infty$ may be verified as in \cite{nolen:ekt04} and \cite{xin:epf92},
and the monotonicity and uniqueness properties follow from the classical parabolic maximum principle, as
described above.

\end{proof}

\end{subsection}

\begin{subsection}{Regularized Equation: Positive Nonlinearity} \label{regularized-kpp}
Now that we have existence of solutions to the regularized problem (\ref{eqn-eps}) for the combustion nonlinearity we wish to
show existence of solutions to (\ref{eqn-eps}) for a KPP nonlinearity $f$.
As in \cite{berestycki:tfc92}, we approximate $f$ by a
sequence of combustion nonlinearities as follows.
Given $\theta \in (0,1)$, let $\chi^{\theta} : \R \rightarrow [0,1]$ be
a smooth cut--off function such that
\br
\left\{ \begin{array}{c}
\displaystyle{ \chi^{\theta}( u ) = 0, \quad \textrm{for} \quad u \leq \theta / 2, } \\
\\
\displaystyle{ \chi^{\theta}( u ) = 1, \quad \textrm{for} \quad u \geq \theta, \quad \textrm{and} } \\
\\
\chi^{\theta} \geq \chi^{\theta'} \quad \textrm{if} \quad \theta \leq \theta' \, .
\end{array} \right.  \label{thetaapprox}
\er
For each $\theta$, the function $f^{\theta} := \chi^{\theta} f $ is a nonlinearity of combustion type,
and the increasing sequence $f^{\theta}$ converges uniformly to $ f $ on $[0,1]$.

For a fixed $\eps > 0$, the results in the preceding section
show that, for each $\theta$, there is a unique classical solution $ ( U^{\theta}, c_{\theta} ) $ of
\begin{equation} \label{kpp-theta}
\left\{ \begin{array}{c}
\displaystyle{ U^{\theta}_{\tau} =  \eps U^{\theta}_{ss} + \left( k \partial_{s} + \nabla_{y} \right)^{2} U^{\theta}
+ b \cdot \left( k \partial_{s} + \nabla_{y} \right) U^{\theta} + c_{\theta} \, U^{\theta}_{s} + f^{\theta}( U^{\theta} ) \, , } \\
\\
\displaystyle{ U^{\theta}(s,y,\tau) \quad 1\textrm{--periodic in} \quad y \quad \textrm{and} \quad  \tau, } \\
\\
\displaystyle{ \lim_{s \to \infty} U^{\theta}(s,y,\tau) = 0 \, , \quad \textrm{and} \quad \lim_{s \to -\infty} U^{\theta}(s,y,\tau) = 1 } \, ,
\\
\end{array} \right.
\end{equation}
which satisfies the normalization condition
\begin{equation} \label{normalize}
\min_{y,\tau}{ \left\{ \, U^{\theta}(0,y,\tau) \, \right\} } = \frac{1}{2} \, .
\end{equation}
Using the parabolic maximum principle, one can show that the wave speeds $c_{\theta}$ increase as $\theta$ approaches $0$,
since the sequence $f^{\theta}$ increases as $\theta$ decreases (\cite{xin:epf92}).
Moreover, these speeds satisfy
\[
c_{\theta} \, \leq \,\overline{c} = \overline{c}_\theta  \, , \quad \textrm{for all} \quad \theta > 0 \, ,
\]
where $\overline{c}_\theta$ was defined above in the proof of Proposition \ref{wave-speed}. Although the constant $\overline{c}_\theta$ depends on $\theta$, the estimates of Section 8 in \cite{berestycki:tfc92} show that the speed of the corresponding one-dimensional traveling waves (\ref{tildepseq}) can be bounded above independently of $\theta$. Hence, the constants $\overline{c}_\theta$ are bounded independently  of $\theta$.  The sequence $c_{\theta}$ therefore converges to $c_{\eps}^{*} > 0$ as $\theta \rightarrow 0$.

Having shown that the speeds $c_{\theta}$ converge as $\theta \rightarrow 0$, we can apply parabolic regularity (as in
the proof of Theorem \ref{existence-eps}) to show that the functions $U^{\theta}$ converge, as $\theta \rightarrow 0$,
to the unique classical solution $U^{\eps}$ of the problem
\begin{equation} \label{kpp-eps}
\left\{ \begin{array}{c}
\displaystyle{ U^{\eps}_{\tau} =  \eps U^{\eps}_{ss} + \left( k \partial_{s} + \nabla_{y} \right)^{2} U^{\eps}
+ b \cdot \left( k \partial_{s} + \nabla_{y} \right) U^{\eps} + c_{\eps}^{*} \, U^{\eps}_{s} + f( U^{\eps} ) \, , } \\
\\
\displaystyle{ U^{\eps}(s,y,\tau) \quad 1\textrm{--periodic in} \quad y \quad \textrm{and} \quad  \tau, } \\
\\
\displaystyle{ \lim_{s \to \infty} U^{\eps}(s,y,\tau) = 0 \, , \quad \textrm{and} \quad \lim_{s \to -\infty} U^{\eps}(s,y,\tau) = 1 } \, ,
\\
\end{array} \right.
\end{equation}
subject to the normalization (\ref{normalize}).  The fact that the limits hold at $s= \pm \infty$ follows from the normalization (\ref{normalize}) and arguments as in \cite{nolen:ekt04} that make use of local uniform convergence.

\end{subsection}

\begin{subsection}{Removal of the Regularization} \label{degenerate-kpp}

The previous section established that, for any given $\eps > 0$, there is a classical solution $(U^{\eps},c_{\eps})$
of problem (\ref{eqn-eps}) for the positive nonlinearity.
To obtain a solution of the original problem (\ref{U}), we want to let the regularization parameter $\eps$ go to $0$.
Therefore, we must prove that the speeds $c_{\eps}$ and the functions $U^{\eps}$ converge in some sense as $\eps \rightarrow 0$.
Proposition \ref{wave-speed} shows that the speeds $c_{\eps}$ are bounded away from $0$ and $+\infty$, while the next result
yields the needed estimates on the functions $U^{\eps}$.

\begin{proposition} \label{w11-est}
The solutions $ U^{\eps} $ of (\ref{eqn-eps}) are uniformly bounded with respect to $\eps$ in $ W^{1,1}_{\loc} $.
\end{proposition}

\begin{proof}

Since $0 < U^{\eps} < 1$ for all $\eps > 0$, $U^{\eps} \in L^{p}_{\loc}$ for any $p \geq 1$.
It follows from the monotonicity of $U^{\eps}$ with respect to $s$ that $U^{\eps}_{s} < 0$ for all $\eps$, and we have
\begin{equation} \label{L1-ds}
\int_{\R \times Q \times (0,1)}{ U^{\eps}_{s} \, ds \, dy \, d \tau } = -1 \quad \textrm{for all} \quad \eps,
\end{equation}
showing that $U^{\eps}_{s} \in L^{1}_{\loc}$.
Multiplying the equation in (\ref{eqn-eps}) by $U^{\eps}$ and integrating (as in the derivation of (\ref{L2-dsdx}) above), we have
\begin{equation} \label{L2-dsdx-e}
\int{ \left| k U^{\eps}_{s} + \nabla_{y} U^{\eps} \right|^{2} \, ds \, dy \, d \tau } =
-\frac{ c_{\eps} }{ 2 } + \int{ U^{\eps} f( U^{\eps} ) \, ds \, dy \, d \tau } \leq
 \frac{c_{\eps}}{2} \leq   \frac{\overline{c}}{2} \, ,
\end{equation}
for the constant $\overline{c}$ from Proposition \ref{wave-speed}.
It follows from (\ref{L2-dsdx-e}) that 

\noindent $ \left| k U^{\eps}_{s} + \nabla_{y} U^{\eps} \right| \in L^{1}_{\loc}$, and we conclude from
(\ref{L1-ds}) that $ \left| \nabla_{y} U^{\eps} \right| \in L^{1}_{\loc} $.  Finally, rearranging equation (\ref{eqn-eps}),
multiplying by $\left( U^{\eps}_{\tau} - c_{\eps} U^{\eps}_{s} \right)$, and integrating, we have
\[
\int{ \left| U^{\eps}_{\tau} - c_{\eps} U^{\eps}_{s} \right|^{2} \, ds \, dy \, d \tau } =
\int{ b \cdot \left( k U^{\eps}_{s} + \nabla_{y} U^{\eps} \right) \left( U^{\eps}_{\tau} - c_{\eps} U^{\eps}_{s} \right) \, ds \, dy \, d \tau } \, .
\]
Applying H\"{o}lder's inequality and using (\ref{L2-dsdx}), we see that $ \left( U^{\eps}_{\tau} - c_{\eps} U^{\eps}_{s} \right) $ belongs to
$L^{2}$, and it follows from (\ref{L1-ds}) that $U^{\eps}_{\tau} \in L^{1}_{\loc}$.  We therefore find that $U^{\eps} \in W^{1,1}_{\loc}$, and
it is clear that, for any compact $E$, the $W^{1,1}(E)$ norm of $U^{\eps}$ can be bounded by a constant independent of $\epsilon > 0$.

\end{proof}

Since the functions $U^{\eps}$ are bounded in $W^{1,1}_{\loc}$, there exists a limit function $U \in L^{1}_{\loc}$
(in fact, $U \in BV_{\loc}$ \cite{evans:mtf92}) such that, after passing to a subsequence,
\[
U^{\eps} \longrightarrow U \quad \textrm{in} \quad L^{1}_{\loc} \quad \textrm{and} \quad U^{\eps} \longrightarrow U \quad \textrm{a.e.} \, .
\]
Of course, the $U^{\eps}$ also converge weakly in $L^{1}_{\loc}$ to $U$; letting $\eps \rightarrow 0$ in equation (\ref{eqn-eps}),
we see that $U$ solves the traveling front equation (\ref{front-eqn}) in the sense of distributions.
%% \begin{equation} \label{front-eqn}
%% U_{\tau} =  \left( k \partial_{s} + \nabla_{y} \right)^{2} U
%% + b \cdot \left( k \partial_{s} + \nabla_{y} \right) U + c \, U_{s} + f( U ) \, .
%% \end{equation}

To analyze this convergence more closely, define the functions $u^{\eps}$ by
\[
u^{\eps}(x,t) := U^{\eps}( k \cdot x - c_{\eps} t, x, t) \, .
\]
The uniform $L^{2}$ estimates on $| k U^{\eps}_{s} + \nabla_{y} U^{\eps} |$ and $ | U^{\eps}_{\tau} - c_{\eps} U^{\eps}_{s} | $
show that the sequence $u^{\eps}$ is compact in $L^{2}_{\loc}$ and thus converges (up to subsequence) to $u \in H^{1}_{\loc}$.
It follows that $u$ is a weak solution of the original semilinear equation (\ref{orig-eqn}), and parabolic regularity results
then imply that $u$ is a classical solution of (\ref{orig-eqn}).  The directional derivatives
\[
U_{\tau} - c U_{s}, \quad
k^{i} U_{s} + U_{y^{i}}, \ i = 1, \ldots, N, \quad \textrm{and} \quad
\left( k \partial_{s} + \nabla_{y} \right)^{2} U
\]
are therefore smooth, but we cannot conclude that $U$ is a classical solution of the
traveling front equation (\ref{front-eqn}).

Now we show that the function $U$ satisfies the desired limits $s \to \pm \infty$:
\be
\lim_{s \to \infty} U(s,y,\tau) = 0 \quad \text{and} \quad \lim_{s \to -\infty} U(s,y,\tau) = 1 .
\ee
The above bounds on $U^\epsilon$ are independent of translation in the $s$ direction, so we are free to normalize the sequence $U^\epsilon$ by
\be
\int_{[0,1] \times Q \times [0,1]} U^\epsilon (s,y,\tau)\,ds\, dy\, d\tau = \frac{1}{2}, \;\;\;\text{for all} \;\epsilon > 0. \label{norm-kpp}
\ee
Since the sequence $U^\epsilon$ converges to $U$ in $L^1_{loc}$, this normalization condition also holds for the limit $U$. Since $U^\epsilon_s < 0$, we conclude that $U$ is nonincreasing in $s$:
\be
U(s,y,\tau) \geq U(s',y,\tau) \;\;\text{if}\;s < s',\;\;\forall (y,\tau) \in Q \times [0,1]. \no
\ee
Therefore, since we also have $0 \leq U \leq 1$, the functions
\be
\psi^+(y,\tau) = \lim_{s \to +\infty} U(s,y,\tau) \quad \text{and} \quad \psi^-(y,\tau) = \lim_{s \to -\infty} U(s,y,\tau) \no
\ee
are well-defined. We claim that $\psi^+(y,\tau) \equiv 0$ and $\psi^-(y,\tau) \equiv 1$. To see this, let $h(y,\tau)$ be any smooth function periodic in $(y,\tau) \in Q \times [0,1]$. For $n=1,2,3 \dots $, define the sequence $\phi_n(s,y,\tau) = \xi_n(s) h(y,\tau)$ where $\xi_n(s) = \xi_0(s - n)$, and $\xi_0$ satisfies
\br
0 \leq \xi_0(s) \leq 1 , \;\;\;\text{for all} \;\;s \in R, \no \\
\xi_0(s) = 0, \;\;\;\text{for}\;\; \abs{s} \geq 1, \label{xi0}
\er
and
\br
\int_R \xi_0(s)\, ds = 1. \no
\er
Note that for each $n$, $\xi_n$ satisfies
\be
\int_R (\xi_n)_s \,ds = 0 = \int_R (\xi_n)_{ss} \,ds.
\ee
Now multiply (\ref{front-eqn}) by $\phi_n$ and integrate over $R \times Q \times [0,T]$. We see that
\[
\int_{R \times Q \times [0,1]} U(s,y,\tau) (\phi^n)_{ss}\,ds\,dy\,d\tau
= \int_{Q \times [0,1]} h(y,\tau) \left( \int_R U (\phi^n)_{ss}\,ds \right) \,dy\,d\tau \]
\[ = \int_{Q \times [0,1]} h(y,\tau) \left( \int_R [\psi^+ + (U - \psi^+)] (\phi^n)_{ss}\,ds \right) \,dy\,d\tau. 
\]
For fixed $(y,\tau)$, we have
\be
\lim_{s \to +\infty} \abs{U(s,y,\tau) - \psi^+(s,y,\tau)} = 0.
\ee
Therefore, for each $(y,\tau)$, we have
\be
\lim_{n \to \infty} \int_R [U - \psi^+] (\phi^n)_{ss}\,ds = 0.
\ee
Also, by definition of $\phi$, for each $(y,\tau)$,
\be
\int_R \psi^+(y,\tau) (\phi^n)_{ss}\,ds = \psi^+(y,\tau) \int_R (\phi^n)_{ss}\,ds  = 0.
\ee
Using the dominated convergence theorem, we now conclude that
\be
\lim_{n \to \infty} \int_{R \times Q \times [0,1]} U (\phi^n)_{ss}\,ds\,dy\,d\tau  = 0.
\ee
Applying a similar argument to each term in the equation, we let $n \to \infty$ and conclude that $\psi^+$ satisfies
\be
\int \psi^+ h_\tau \,dy\,d\tau  - \int \psi^+ \Delta_y h \,dy\,d\tau - \int \psi^+ b \cdot \nabla_y h \,dy\,d\tau = \int f(\psi^+) h \,dy\,d\tau .
\ee
Since $h$ was arbitrary, we conclude that $\psi^+$ is a distributional solution of the periodic-parabolic equation
\be
\psi_\tau - \Delta_y \psi + b \cdot \nabla_y \psi = f(\psi) \label{psi-eqn}
\ee
Because of the bounds on $U$, we know that $0 \leq \psi^+ \leq 1$, so that right hand side $f(\psi^+) \in L^p(Q \times [0,1])$, for any $p > 0$. By parabolic regularity estimates, it now follows that $\psi^+ \in W^{2,p}$ for any $p>0$. Then, as in the analysis following (\ref{wtwopest}), we find that $\psi^+$ is a classical solution of (\ref{psi-eqn}), assuming $b$ is sufficiently smooth. The maximum principle implies that $\psi^+$ must be constant. Hence, $f(\psi^+) \equiv 0$, so that $\psi^+ \equiv 0$ or $\psi^+ \equiv 1$. Because of the normalization (\ref{norm-kpp}) and because $U$ is nonincreasing, $\psi^+ \leq \frac{1}{2}$. Therefore, $\psi^+ \equiv 0$. A similary argument with $\xi_n(s) = \xi_0(s + n)$ implies that $\psi^- \equiv 1$.

We have now shown that $U$ satisfies the desired limits, concluding the proof of Theorem \ref{existence-kpp}.

\end{subsection}

\end{section}

\begin{section}{Propagation Speed for KPP Nonlinearity} \label{section-existence-kpp}
\setcounter{equation}{0}
In this section, we make the additional assumption that $f$ is a KPP
nonlinearity satisfying $f(u) \leq u f'(0)$.  Under this assumption, we find that the speed $c^*$ of the traveling wave obtained in the preceding construction is minimal and can be characterized by a variational formula which extends the well-known variational formula for the KPP minimal speed when $b$ is time-independent. Moreover, we can describe propagation of front-like solutions at any speed higher than $c^*$. First, we characterize the speed $c^*$.

\begin{subsection}{Characterization of the KPP Minimal Speed $c^*$}
To introduce the variational characterization, we define $\mu(\lambda)$ to be the principal eigenvalue of the periodic--parabolic operator (\cite{hess:ppb91})
\begin{equation} \label{eigen-op}
L^{\lambda} \Phi :=
\Delta_{x} \Phi + \left( b - 2 \lambda k \right) \cdot \nabla_{x}\Phi +
\left( \lambda^{2} - \lambda b \cdot k + f'(0) \right) \Phi - \Phi_{t} \, .
\end{equation}
defined on spatially--temporally periodic functions $\Phi(x,t)$. The associated eigenfunction $\Phi = \Phi_\lambda > 0$ is unique up to multiplication by a constant.
We have the following characterization of the speed $c^*$ obtained in the preceding construction:
\begin{theorem}
The speed $c^*$ of the traveling wave obtained in the preceding sections satisfies
\begin{equation} \label{min-speed}
c^{*} =  \inf_{\lambda > 0}{ \frac{ \mu( \lambda ) }{ \lambda } }
\end{equation}
Moreover, this speed is minimal in the sense that if  $(U,c)$ is any traveling front solution of (\ref{front-eqn}), we must have $c \geq c^*$.
\end{theorem}
\begin{proof}
In the special case that $b$ is a time-dependent shear flow, this characterization was proven in \cite{nolen:ekt04}. The proof there generalized ideas of \cite{BH1} and relied heavily on regularity properties of the traveling wave equation. In the present case, the traveling wave equation is degenerate. Nevertheless, for the regularized equation, the arguments of \cite{nolen:ekt04} can be extended in a straightforward way to prove that
\begin{equation} \label{min-speed-eps}
c_{\eps}^{*} =  \inf_{\lambda > 0}{ \frac{ \mu^{\eps}( \lambda ) }{ \lambda } } \, ,
\end{equation}
where $\mu^{\eps}( \lambda ) $ denotes the principal eigenvalue of the periodic--parabolic operator
\begin{equation} \label{eigen-op-eps}
L^{\epsilon,\lambda} \Phi :=
\Delta_{x} \Phi + \left( b - 2 \lambda k \right) \cdot \nabla_{x}\Phi +
\left( (1 + \epsilon)\lambda^{2} - \lambda b \cdot k + f'(0) \right) \Phi - \Phi_{t} \, ,
\end{equation}
defined on spatially--temporally periodic functions $\Phi(x,t)$.
Note that
\begin{equation} \label{mu-eps}
\mu^{\eps}( \lambda ) = \mu( \lambda ) + \eps \lambda^{2}  .
\end{equation}
so that
\begin{equation}
\inf_{\lambda > 0} \frac{\mu^\epsilon(\lambda)}{\lambda} = \inf_{\lambda > 0} \left( \frac{\mu(\lambda)}{\lambda} +\epsilon \lambda \right )
\end{equation}
Using  the properties of $\mu(\lambda)$ described in the appendix, we see that the speeds $c_{\eps}^{*}$ converge to
\begin{equation} \label{min-speed1}
c^{*} =  \inf_{\lambda > 0}{ \frac{ \mu( \lambda ) }{ \lambda } }
\end{equation}
as $\eps \rightarrow 0$.  We will show below that this limit $c_{*}$ is indeed the minimal speed at which
a traveling wave solution of (\ref{orig-eqn}) may propagate.
\end{proof}

\end{subsection}

\begin{subsection}{Generalized Propagation and Minimality of $c^*$}
In this section we show that the speed $c^*$ defined above is indeed the minimal speed of propagation and that solutions of the corresponding initial value problem can propagate at speeds higher than $c^*$.  To do so, we generalize some results found in \cite{AW} and \cite{MR} that define the asymptotic propagation speed in the KPP case, without relying on the existence of smooth traveling waves. We consider propagation in a given direction $k$, a unit vector in $R^N$, arising from solutions to the initial value problem
\br
L u + f(u) = \Delta_x u - u_t + b(x,t) \cdot \nabla u + f(u) = 0 \label{e5} \\
u(x,0) = u_0(x). \no
\er
where $f$ is the KPP nonlinearity satisfying $f(u) \leq f'(0)u$ for $u \in (0,1)$. For a given $r \in R$, we use $\Sigma^-_r \subset R^N$ to denote the open half-space:
\be
\Sigma^-_r = \left\{ x \in R^N \;| \;\; k\cdot x < r\right\}. \no
\ee
Similarly we define $\Sigma^+_r$ by $(k\cdot x) > r$, i.e. the interior of the complement of $\Sigma^-_r$.  Also, we define for $r > 0$,
\be
\Sigma^0_r = \left\{ x \in R^N \;| \;\; \abs{k\cdot x} < r\right\}. \no
\ee
Unless otherwise stated, we will consider initial data $u_0$ satisfying $0 \leq u_0 \leq 1$ and the limits
\be
\lim_{r \to -\infty} \inf_{\Sigma^-_r} u_0(x) = 1 \;\;\; \text{and} \;\;\; \lim_{r \to +\infty} \sup_{\Sigma^+_r} u_0(x) = 0 \label{e6}
\ee
That is, $u_0$ will converge to limits $0$ and $1$ uniformly in the directions $k$ and $-k$, respectively. In this case we say $u_0$ is wave-like in the direction $k$. The main tool for this analysis is the eigenvalue problem associated to propagation in the direction $k$:
\be
L_\lambda \phi = \Delta \phi - \phi_t + (b - 2\lambda k)\cdot \nabla \phi + (\lambda^2 - \lambda(b\cdot k) + f'(0) )\phi = \mu(\lambda) \phi \label{e1}
\ee
with $\phi(x,t) = \phi_\lambda > 0$ being periodic in both $x$ and $t$. This is obtained by substituting a solution of the form $e^{-\lambda(k\cdot x - ct) } \phi(x,t)$ into the linearized version of equation (\ref{e5}) and taking $\mu(\lambda) = \lambda c$. Define $c^*$ by
\be
c^* = \inf_{\lambda > 0} \frac{\mu(\lambda)}{\lambda} > 0. \label{e3}
\ee
Let $\lambda^* > 0$ be the point where this infimum is attained.  As described in the appendix, for each $c > c*$ there exists a unique $\lambda_c \in (0, \lambda^*)$ such that $\mu(\lambda_c) = c \lambda_c$.

Note that $c^*$ and $\lambda^*$ depend on the direction vector $k$. If we consider propagation in the opposite direction, $k \to -k$, then (\ref{e3}) gives us another $c^{**}$ and $\lambda^{**}$ corresponding to $-k$, and it is not necessarily true that $c^{**} = c^*$.

\begin{definition}[Speed of Propagation] \label{sop}
We say that the solution to the initial value problem (\ref{e5}) propagates in direction $k$ with speed $c$ if for any $c' > c$ and any $r \in R $,
\be
\lim_{t \to \infty} \sup_{x \in \Sigma^+_r} u(x + c'tk,t) = 0, \no
\ee
and for any $c'' < c$ and any $r \in R$,
\be
\lim_{t \to \infty} \inf_{x \in \Sigma^-_r} u(x + c''tk,t) = 1. \no
\ee
\end{definition}

The following theorems and their corollaries allow us to characterize the asymptotic speed in terms of the decay rate of the initial data.  Here are the main results:

\begin{theorem} \label{th:sop1}
Let $u_0$ satisfy (\ref{e6}) and $u_0(x) \leq C e^{-\lambda_c k\cdot x}$ for some constant and $c \geq c^*$.  Then for any $c' > c$ and any $r \in R$ we have
\be
\lim_{t \to \infty} \sup_{x \in \Sigma^+_r} u(x + c'tk,t) = 0.
\ee
That is, the solution $u$ propogates with speed $\leq c$.
\end{theorem}

\begin{theorem} \label{th:sop2}
If there exists a constant $m>0$ and an interval $[a_1, a_2] \subset R$ such that
\br
0 \leq u_0(x) \leq 1 \;\;\text{for all}\;\; x\in R^N, \no \\
u_0(x) > m > 0 \;\;\text{wherever}\;\; (k\cdot x) \in [a_1,a_2],
\er
then
\be
\lim_{t \to \infty} \inf_{\Sigma^0_r} u(x + c'kt, t) = 1
\ee
for any $c' \in (-c^{**}, c^{*})$ and $r \in R$.
\end{theorem}

\begin{theorem} \label{th:sop3}
Let $c > c^*$ and let $u_0$ satisfy (\ref{e6}) and $u_0(x) \geq C_1 e^{-\lambda_c k\cdot x}$ for $(k \cdot x) > C_2$, for some constants $C_1, C_2$. Then for any $c' < c$, we have
\be
\lim_{t \to \infty} \inf_{x \in \Sigma^-_r} u(x + c'tk,t) = 1.
\ee
\end{theorem}

\begin{corollary} \label{cor:sop1}
For initial data satisfying (\ref{e6}),
\be
\lim_{t \to \infty} \inf_{x \in \Sigma^-_r} u(x + c'tk,t) = 1
\ee
for any $c' < c^*$.
\end{corollary}

\begin{corollary} \label{cor:sop2}
If $u_0$ satisfies (\ref{e6}) and $u_0(x) \leq Ce^{-\lambda^* k\cdot x}$ for some constant $C>0$, then the asymptotic speed of propagation in direction $k$ must be $c^*$.
\end{corollary}

Recall that $-c^{**}$ corresponds to propagation in the $-k$ direction at speed $\abs{-c^{**}}$.  The point is that an initial disturbance must spread at a rate at least as fast as the minimal speeds in the corresponding direction. In what follows, we rely heavily on the eigenvalue problem (\ref{e1}) to construct global sub- and super-solutions to the initial value problem (\ref{e5}).

\nit {\bf Proof of Theorem \ref{th:sop1}:}
For any $c_1 \in (c, c')$, we have $c_1 \lambda_c > \mu(\lambda_c)$, by definition of $\lambda_c$.  Then define $\psi(x,t) = M e^{-\lambda_c (k \cdot x - c_1 t)} \phi(x,t)$, where $\phi = \phi_{\lambda_c} > 0$ is defined by (\ref{e1}), and the constant $M$ is chosen so that $u_0(x) < \psi(x,0)$ for all $x \in R^N$. Then $\psi$ is a supersolution:
\br
L \psi + f'(0) \psi < 0 \no \\
\psi(x,0) > u_0(x). \no
\er
Since $f(\psi) \leq \psi f'(0)$, it follows from the maximum principle that $\psi(x,t) > u(x,t)$ for all $ x \in R^N$, $t \geq 0$. Therefore,
\br
\lim_{t \to \infty} \sup_{x \in \Sigma^+_r} u(x + c'tk,t) &\leq& \lim_{t \to \infty} \sup_{x \in \Sigma^+_r} \psi(x + c'tk,t) \no \\
& = &  \lim_{t \to \infty} \sup_{x \in \Sigma^+_r} M e^{-\lambda_c (k \cdot x + (c' - c_1) t)} \phi(x + c'tk,t) \no \\
& \leq & M \norm{\phi}_\infty e^{-\lambda_c r} \lim_{t \to \infty}  e^{-\lambda_c (c' - c_1)t} = 0,
\er
since $c' > c_1$. This proves Theorem \ref{th:sop1}.

The proof of Theorem \ref{th:sop2}  extends ideas from Theorem 3.4 of \cite{MR}. First, we state the following lemma, which is a straightforward extension of Lemma 3.5 in \cite{MR}:
\begin{lemma} \label{lem:sop1}
There exists a $c' < c^*$ such that for all $c \in (c',c^*)$, there exists $\lambda \in \mathcal{C} \setminus \mathcal{R}$ and a function $\phi = \phi(x,t) \in C^{2,1}(T^{n+1},\mathcal{C})$ with $Re(\phi) > 0$ and solving
\be
\Delta \phi - \phi_t + (b - 2\lambda k)\cdot \nabla \phi + (\lambda^2 -\lambda c - \lambda(b\cdot k) + f'(0) )\phi = 0 .\label{e2}
\ee
\end{lemma}
The important point of the lemma is that we can choose $\lambda$ to be strictly complex and $\phi$ to have strictly positive real part.  This will allow us to construct global subsolutions that have no tail (or sufficiently small tail) and can fit beneath $u_0(x)$.  We omit the proof of this lemma since it is identical to that of Lemma 3.5 of \cite{MR}, except that the operator is somewhat different because of the time-dependence and more general flow structure.

\nit {\bf Proof of Theorem \ref{th:sop2}:}

{\bf Step 1:} If we replace $f'(0)$ by $f'(0) - \delta$ in (\ref{e1}), then from (\ref{e3}) we obtain a new minimal speed $c^*_\delta$ instead of $c^*$ with $c^*_\delta \nearrow c^*$ as $\delta \to 0$. Also the lemma would hold with $f'(0)$ replaced by $f'(0) - \delta$, so that for $\delta$ sufficiently small and for $c_R < c^*_\delta < c^*$ sufficiently close to $c^*_\delta$, we apply the lemma to obtain a function $\phi$ such that
\be
\Delta \phi - \phi_t + (b - 2\lambda k)\cdot \nabla \phi + (\lambda^2 -\lambda c_R - \lambda(b\cdot k) + f'(0) - \delta )\phi = 0 ,\no
\ee
where $Re(\phi) > 0$ and $\lambda \in \mathcal{C} \setminus \mathcal{R}$ depends on $c_R$. Now let
\be
\psi = Re(e^{-\lambda (k \cdot x - c_R t)} \phi (x,t)). \no
\ee
Then $\psi$ solves
\be
 \psi_t - \Delta \psi - b \cdot \nabla \psi = (f'(0) - \delta)\psi, \no
\ee
or in the moving frame $\hat \psi(x,t) = \psi(x + c_R kt,t)$:
\be
\hat \psi_t - \Delta \hat \psi -  (b + c_R k) \cdot \nabla \hat \psi = (f'(0) - \delta)\hat \psi. \no
\ee
Now
\be
\psi(x,t) = e^{-\lambda_r (k\cdot x - c_R t)} \left[ \phi_r \cos(\lambda_i (k\cdot x - c_R t)) + \phi_i \sin(\lambda_i(k\cdot x - c_R t))  \right] \label{e4}
\ee
where $\lambda_r, \lambda_i, \phi_r, \phi_i$ denote the real and imaginary parts of $\lambda$ and $\phi$. If $(k\cdot x - c_R t) = \frac{ 2 n \pi}{\lambda_i}$ for $n \in \mathcal{Z}$, then $\psi > 0$, since $\phi_r > 0$. Similarly, if $(k\cdot x - c_R t) = \frac{ (2n+1) \pi}{\lambda_i}$, then $\psi < 0$. It follows from (\ref{e4}), i.e. from the complexity of $\lambda$, that for a given $\eta > 0$, there exists an unbounded set $D \subset R^{N+1}$ and an interval $(a_1, a_2) \subset R$ such that
\br
0 < \psi(x,t) < \eta \;\;\; \text{for}\;(x,t) \in D , \no \\
\psi(x,t) = 0 \;\;\; \text{for}\;(x,t) \in \partial D , \no \\
D \subset \left\{ (x,t) \;|\;\; (k\cdot x - c_R t) \in (a_1,a_2)\right\}. \no
\er
Furthermore, we may choose some $a_3 \in (a_1, a_2)$ such that $a_3 = \frac{2n\pi}{\lambda_i}$ for some $n \in \mathcal{Z}$. Now we define
\be
\psi_R(x,t) = \left\{ \begin{array}{c}
\psi(x,t) \;\; \text{if}\;(x,t) \in D \\
\\
0 \;\; \text{otherwise} .\\
\end{array} \right.
\ee
So, $0 \leq \psi_R < \eta$ with support $D$ bounded by the two hyperplanes $(k \cdot x - c_R t) = a_1$ and $(k \cdot x - c_R t) = a_2$. Also, $\psi_R(x,t) >0$ whenever $(k\cdot x - c_R t) = a_3$. If we choose $\eta > 0$ such that $f(s) > s(f'(0)-\delta)$ for $s \in (0,\eta]$, then $\psi_R$ satisfies
\be
(\psi_R)_t - \Delta \psi_R - b\cdot \nabla \psi_R \leq f(\psi_R) \no
\ee
for all $(x,t) \in D$. By translating $\psi_R$ in the variable $t$, we may assume that $\psi_R$ has support in the $(a_1,a_2)$ defined in the statement of the theorem, so that by reducing $\eta$, we may assume $\psi_R(x,0) \leq u_0(x)$ for all $x \in R^N$.

Now $\psi_R$ is the desired subsolution. The function $u(x,t)$ satisfies
\br
u_t - \Delta u - b \cdot \nabla u = f(u) \no \\
u(x,0) =  u_0(x) \geq \psi_R(x,0) .\no
\er
By maximum principle $u > 0$ for $t>0$, and $u \geq \psi_R $ for all $t \geq 0$.

Suppose we repeat the preceeding analysis after switching $k \to -k$. In this case, we have a new minimal speed $c^{**}_\delta$ corresponding to waves in the direction $-k$. Furthermore, we can create a bump function $\psi_L$ (analogous to $\psi_R$) that is a subsolution propagating in the $-k$ direction with speed $c_L$ arbitrarily close to $c^{**}_\delta$ and $-c^{**}_\delta < -c_L < c_R < c^*_\delta$. It follows that if $u(x,t)$ satisfies
\br
u_t - \Delta u - b \cdot \nabla u = f(u) \no \\
u_0(x) \geq \psi_R(x,0) \;\; \text{and} \;\;\; u_0(x) \geq \psi_L(x,0) \, ,\no
\er
then $u \geq \psi_R$ and $u \geq \psi_L$ for all $t \geq 0$. In particular, there exists $\epsilon > 0, \Lambda > 0$ and constants $a^0_L, a^0_R$ such that
\br
u(x,t) \geq \psi_L(x,t) > \epsilon > 0 \;\; \text{if}\;\; \abs{k \cdot x - a_L(t)} < \Lambda \no \\
u(x,t) \geq \psi_R(x,t) > \epsilon > 0 \;\; \text{if}\;\; \abs{k \cdot x - a_R(t)} < \Lambda, \label{e8}
\er
where $a_L(t) = a^0_L - c_L t$ and $a_R(t) = a^0_R + c_R t$. The constant $\epsilon$ is determined by $\eta$ and $\inf_{(x,t)} \phi_{r,R}$, and $\inf_{(x,t)} \phi_{r,L}$, where $\phi_{r,R}$ and $\phi_{r,L}$ are the real parts of $\phi$ defined in the construction of $\psi_R$ and $\psi_L$, as in (\ref{e4}).

{\bf Step 2:}
Now with $a_L(t), a_R(t)$ defined as above, we define for $\tau \geq 0$ the set $\Omega_\tau \subset R^{N+1}$ by
\be
\Omega_{\tau} = \left\{ (x,t) | \; k \cdot x \in [a_L(t), a_R(t)], t \geq \tau \right\}. \no
\ee
and we show that
\be
\lim_{\tau \to \infty} \inf_{\Omega_\tau} u(x,t) \geq \epsilon. \label{e9}
\ee
In other words, $u$ stays bounded uniformly away from zero in a region spreading at rates $c_L$ and $c_R$ in the directions $-k$ and $k$, respectively.

Let $\hat f \in C^1([0,1])$ be a nonlinear function satisfying $\hat f(u) \leq f(u)$ for all $u \in [0,1]$ and
\br
\hat f(u) = 0 \;\;\text{for}\;\; u \in \{ 0 \} \cup [\epsilon, 1], \no \\
\hat f(u) > 0 \;\; \text{for}\;\; u \in (0,\epsilon). \no
\er
Now let $t_0 > 0$ such that $a_L(t) < a_R(t)$ for $t \geq t_0$ and take
\be
g_0 = \frac{1}{2} \inf_x \left\{ u(x,t_0) |\; (k\cdot x) \in [a_L(t_0), a_R(t_0)] \right \} \no
\ee
Let the function $g(t)$ solve $g'(t) = \hat f(g(t))$ for $t \geq t_0$ and $g(t_0) = g_0$. Then in the region $\Omega_{t_0}$, we have
\br
g_t - \Delta g - b \cdot \nabla g = \hat f(g) < f(g) \no \\
u_t - \Delta u - b \cdot \nabla u = f(u) \no
\er
with $u(x,t) > g(t)$ in the parabolic boundary $\partial \Omega_{t_o}$, since by (\ref{e8}) $u > \epsilon$ on the lateral boundary of $\Omega_{t_o}$. Maximum principle now implies that $u(x,t) > g(t)$ for all $(x,t) \in \Omega_{t_o}$. Since $g(t) \to \epsilon$, this implies the (\ref{e9}).

{\bf Step 3:}
In this step we improve the result of Step 2 by showing that for any $q \in (0,1/2)$ arbitrarily small, then for some constant $d = d(q)$ to be determined,
\be
\lim_{\tau \to \infty} \inf_{\Omega'_\tau} u(x,t) \geq 1 - q \label{clstep3}
\ee
where
\be
\Omega'_{\tau} = \left\{ (x,t) | \; k \cdot x \in [a_L(t) + d, a_R(t) - d], t \geq \tau \right\}. \no
\ee
Inspired by the technique of \cite{ES1} (see p. 158), we let $\theta \in (0,1)$ and define the linear function
\be
\bar f(u) = s(1- \theta - u) \no
\ee
choosing the constant $s > 0$ so that $f(u) > \bar f(u)$ whenever $u \geq \epsilon$. Thus by the preceeding step, there is $t_1>0$ sufficiently large so that
\be
u_t - \Delta u - b \cdot \nabla u > \bar f(u) \;\;\; \text{for all}\;\; (x,t) \in \Omega_{t_1}. \no
\ee
For constants $\gamma, \kappa \in (0,1)$ to be determined, define
\br
\bar a_L(t) = a_L(t) + \left(\frac{\kappa}{\gamma}\right)^{\frac{1}{4}}, \no \\
\bar a_R(t) = a_R(t) - \left(\frac{\kappa}{\gamma}\right)^{\frac{1}{4}}. \no
\er
and the function
\br
\rho(x,t) = \left\{ \begin{array}{c}
\kappa, \;\;\; \text{if}\; (k\cdot x) \in [\bar a_L(t), \bar a_R(t)]  \\
- \gamma (k\cdot x - \bar a_R(t))^4 + \kappa, \;\;\; \text{if}\; (k\cdot x) > \bar a_R (t) \\
- \gamma (k\cdot x - \bar a_L(t))^4 + \kappa, \;\;\; \text{if}\; (k\cdot x) < \bar a_L(t) \\
\end{array} \right.
\er
Note that the functions $\bar a_L(t)$ and $\bar a_R(t)$ were defined so that $\rho(x,t) \geq 0$ if $(x,t) \in \Omega_t$ and $ \rho < 0$ otherwise. Choose $\sigma \in (0,s)$ and let
\be
g(t) = 1 - \theta - e^{-\sigma (t-t_1)} \no
\ee
Now we compare with $u$ the function defined by
\be
z(x,t) = g(t) \rho(x,t) \no
\ee
We compute:
\br
& & z_t - \Delta z - b \cdot \nabla z - \bar f(z) \no \\
& =  & g'(t) \rho + g(t) \rho_t - g(t) (\Delta \rho + b \cdot \nabla \rho) - s(1 - \theta - g \rho) \no \\
&=& (\sigma - s) e^{- \sigma (t - t_1)}\rho + g(t)(\rho_t -  \Delta \rho - b \cdot \nabla \rho) + s(1-\theta)(\rho - 1) \no \\
&\leq & (1-\theta) (\rho_t -  \Delta \rho - b \cdot \nabla \rho) + s(1-\theta)(\kappa - 1)
\er
 A straightforward computation shows that if $(x,t) \in \Omega_t$, where $\rho \geq 0$, then there is a constant $C>0$ such that
\be
\abs{\rho_t} \leq C \gamma^{\frac{1}{4}}, \;\;\; \abs{\nabla \rho} \leq C \gamma^{\frac{1}{4}}, \;\;\; \abs{\Delta \rho} \leq C \gamma^{\frac{1}{4}} \no
\ee
Thus choosing $\gamma$ sufficiently small, we can make $ (1-\theta) (\rho_t -  \Delta \rho - b \cdot \nabla \rho)$ arbitrarily small, so that $z_t - \Delta z - b \cdot \nabla z - \bar f(z) < 0$ since $\kappa < 1$. Therefore, in the region $\Omega_{t_1}$,
\br
z_t - \Delta z - b \cdot \nabla z < \bar f(z) \no \\
u_t - \Delta u - b \cdot \nabla u > \bar f(u)
\er
with $u(x,t)>z(x,t)$ on the parabolic boundary $\partial \Omega_{t_1} $.  Maximum principle implies that $u(x,t) > z(x,t)$ for all $(x,t) \in \Omega_{t_1}$. Therefore,
\be
\lim_{\tau \to \infty} \inf_{\Omega'_\tau} u(x,t) \geq \lim_{\tau \to \infty} \inf_{\Omega'_\tau} z(x,t) = \kappa (1-\theta)
\ee
where
\be
\Omega'_{\tau} = \left\{ (x,t) | \; k \cdot x \in [a_L(t) + \left(\frac{\kappa}{\gamma}\right)^{\frac{1}{4}}, a_R(t) - \left(\frac{\kappa}{\gamma}\right)^{\frac{1}{4}}], t \geq \tau \right\}.
\ee
Since $\theta$ and $\kappa$ can be chosen arbitrarily close to $0$ and $1$, respectively, we can make $\kappa (1 - \theta) > 1-q$, proving the claim (\ref{clstep3}).

{\bf Step 4:}
Now we complete the proof. Recall that $a_L(t) = a^0_L - c_L t$ and $a_R(t) = a^0_R + c_R t$. Therefore, the result (\ref{clstep3}) of preceding step implies that for any $c' \in (-c_L, c_R)$ we have
\be
\lim_{t \to \infty} \inf_{\Sigma^0_r} u(x + c'kt, t) = 1 \label{st4res}
\ee
By definition of $c_L$ and $c_R$, we can choose $c_L$ and $c_R$ arbitrarily close to $c^{**}_\delta$ and $c^*_\delta$, respectively. Furthermore, we can choose $\delta$ sufficiently small so that $c^{**}_\delta$ and $c^*_\delta$ are arbitrarily close to $c^{**}$ and $c^*$, respectively. Hence the result (\ref{st4res}) applies to any $c' \in (-c^{**}, c^*)$.
This completes the proof of Theorem \ref{th:sop2}.

\nit {\bf Proof of Theorem \ref{th:sop3}:}
Without loss of generality, assume $C_2 = 0$. Thus, $u_0(x) \geq C e^{-\lambda_c k\cdot x}$ for $(k \cdot x) > 0$.  We now construct a subsolution of the form
\be
\psi(x,t) = d_1 e^{-\lambda_1 (k \cdot x - c_2 t)} \phi_{\lambda_1}(x,t) - d_2 e^{-\lambda_2 (k \cdot x - c_2 t)} \phi_{\lambda'} (x,t)
\ee
for some $\lambda_1, \lambda_2$ satisfying $\lambda_c < \lambda_1 < \lambda_2$. As described in the appendix, if we replace $f'(0)$ with $f'(0) - \delta$ in the eigenvalue problem, minimal speed $c^*$ becomes $c^*_\delta$, with $c^*_\delta < c^*$, $c^*_\delta \to c^*$ as $\delta \to 0$. Let $\mu_\delta(\lambda)$ denonte the eigenvalue corresponding to the problem with $f'(0) - \delta$.  Then, we can pick $\delta$ small so that there exists $\lambda_1, \lambda_2$ and $c_1, c_2$ such that
\be
c_2 = \frac{\mu_\delta(\lambda_2)}{\lambda_2} <  \frac{\mu_\delta(\lambda_1)}{\lambda_1} = c_1 \label{e7}
\ee
and $ c_2 < c_1 < c$ and $\lambda_c < \lambda_1 < \lambda_2$ (see appendix for properties of $\mu_\delta$). Now define the function
\be
\psi(x,t) = d_1 e^{-\lambda_1 (k \cdot x - c_2 t)} \phi_1(x,t) - d_2 e^{-\lambda_2 (k \cdot x - c_2 t)} \phi_2 (x,t)
\ee
where $\phi_1, \phi_2 > 0$ are the eigenfunctions defined by (\ref{e1}) with $f'(0) - \delta$, corresponding to $\lambda_1, \lambda_2$.
Given an $\eta > 0$, we can pick $d_1, d_2 > 0$ such that $\psi(x,0) < u_0(x)$ for all $x \in R^N$ and $\sup_{(x,t)} \psi(x,t) < \eta$. This is possible since $u_0(x) \geq C e^{-\lambda_c k\cdot x}$ for $(k \cdot x) > 0$ and $u_0$ satisfies (\ref{e6}) and $\lambda_c < \lambda_1 < \lambda_2$.  Also, for $r$ sufficiently large,
\be
\psi(x,0) > 0 \;\; \text{for} \;\; x \in \Sigma^+_r . \no
\ee
By (\ref{e7}), $ \lambda_1 c_2  \phi_1 < \mu_\delta(\lambda_1) \phi_1$. As a result, $\psi$ is a subsolution to the linearized problem:
\br
\psi_t - \Delta \psi - b \cdot \nabla \psi < (f'(0)- \delta) \psi  \no \\
\psi(x,0) < u_0(x) .\no
\er
Now, since $f(s) > s (f'(0) - \delta)$ for $s \in (0,\eta)$ and since $\psi < \eta$ for all $(x,t)$,
\br
\psi_t - \Delta \psi - b \cdot \nabla \psi <  f(\psi)  \no \\
\psi(x,0) < u_0(x). \no
\er
By maximum principle, $u \geq \psi$ for all $(x,t)$. As in (\ref{e8}), we find conclude that there is an $\epsilon > 0$ and $\Lambda > 0$ and a constant $a_0$ such that
\be
u(x,t) \geq \psi(x,t) > \epsilon > 0 \;\; \text{if}\;\; \abs{k \cdot x - a_R(t)} < \Lambda\no \\
\ee
where $a_R(t) = a_0 + c_2 t$.

Now using the assumption that
\be
\lim_{r \to -\infty} \inf_{\Sigma^-_r} u_0(x) = 1, \no
\ee
we proceed as in Step 2 of Theorem \ref{th:sop2} to show that
\be
\lim_{\tau \to \infty} \inf_{\Omega_\tau} u(x,t) \geq \epsilon \no
\ee
on the set
\be
\Omega_{\tau} = \left\{ (x,t) | \; k \cdot x \in (-\infty, a_R(t)], t \geq \tau \right\}. 
\ee
Then it is easy to extend the arguments of Step 3 and Step 4 in the proof of Theorem \ref{th:sop2} to conclude the proof that
\be
\lim_{t \to \infty} \inf_{x \in \Sigma^-_r} u(x + c'tk,t) = 1.
\ee

\nit {\bf Proof of Corollary \ref{cor:sop1}:}  Suppose the initial data is wave-like, i.e. satisfies (\ref{e6}). Without loss of generality, we may assume that $u_0(x) > 0$ if $(k \cdot x) \leq 0$. Then let $v(s)\geq 0$ be a smooth function with compact support in the interval $s \in (0,1)$, with $\norm{v}_\infty = 1$. Then consider the function $\psi_0(x)$ defined by
\be
\psi_0(x) = \sum_{j=1}^\infty a_j v(k\cdot x + j)
\ee
where the coefficients $\left \{ a_j \right\}$ are chosen such that $\psi(x) < u_0(x)$.  Applying the argument from Theorem \ref{th:sop2}, the solution to the initial value problem
\br
L \psi + f(\psi) = 0 \no \\
\psi(x,0) = \psi_0
\er
must satisfy
\be
\lim_{t \to \infty} \inf_{x \in \Sigma^-_r} \psi(x + c'tk,t) = 1
\ee
for any $c' < c^*$. By maximum principle, the result holds for $u(x,t)$ since $u_0(x) \geq \psi_0$.

\nit {\bf Proof of Corollary \ref{cor:sop2}:}
This follows from the preceeding corollary and from Theorem 1.1, with $\lambda_c = \lambda^*$.

\end{subsection}
\end{section}

\section{Appendix: Properties of the Eigenvalue Problem}
The eigenvalue problem associated with propagation in the direction $k$ is:
\be
L_\lambda \phi = \Delta \phi - \phi_t + (b - 2\lambda k)\cdot \nabla \phi + (\lambda^2 - \lambda(b\cdot k) + f'(0) )\phi = \mu(\lambda) \phi \label{ev6}
\ee
with $\phi(x,t) = \phi_\lambda > 0$ being periodic in both $x$ and $t$. In this section we describe the following properties of the principal eigenvalue $\mu(\lambda)$ as a function of $\lambda$:

\begin{lemma} \label{lem:ev1}
The map $\lambda \mapsto \mu(\lambda)$ is strictly convex in $\lambda$.
\end{lemma}
\begin{proposition} \label{prop:ev1}
There is a point $\lambda^* \in (0,+\infty)$ at which the curve $\lambda \mapsto \frac{\mu(\lambda)}{\lambda}$ achieves a unique global minimum (over $(0,+\infty)$) which is also the only local minimum. Hence, $\frac{\mu(\lambda)}{\lambda}$ is strictly decreasing over interval $(0,\lambda^*)$ and increasing over interval $(\lambda^*, +\infty)$.
\end{proposition}
{\bf Proof of Proposition \ref{prop:ev1}}
Integrate the equation (\ref{e1}) over $D = Q \times [0,T]$ and divide by $\int_D \phi \;dx\,dt > 0$.  Periodicity and the fact that $\nabla \cdot b = 0$ imply that
\be
\mu(\lambda) = \lambda^2 + f'(0) + \lambda \frac{\int_D (b \cdot k) \phi}{ \int_D \phi}.
\ee
Since $\phi > 0$, the last term is bounded by $\norm{b}_\infty$, so that $\mu(\lambda)$ grows quadratically as $\lambda \to \pm \infty$. Letting $\lambda=0$, maximum principle implies that $\mu(0) = f'(0)>0$ with eigenfunction being a positive constant. Then the proposition follows immediately from the fact that $\mu(\lambda)$ is convex in $\lambda$.

The convexity of $\mu(\lambda)$ follows from a minor extension of the argument in \cite{BH1} (see Section 5.2 therein), but for convenience and since the eigenvalue problem is nonstandard, we include the details here. First, we make use of the following characterization of the principal eigenvalue, which holds even though the operator is not self-adjoint. Let operator $L$ be defined by
\be
L \phi = a_{ij} \phi_{x_i x_j} - \phi_t + b_i \phi_{x_i} + c \phi
\ee
where $a_{ij} = (a_{ij}(x,t))$ is uniformly elliptic, $b_i = b_i(x,t)$ and $c = c(x,t)$, with summation notation implied. All coefficients are periodic with respect to both $x$ and $t$. Denote by $\mu$ the principle eigenvalue of $L$ over the space $E =  C_P^2(Q \times [0,T])$ of $C^2$ functions periodic in $D = Q \times [0,T]$.

\begin{theorem} \label{th:ev1}
The following variational formula holds for the principal eigenvalue $\mu$:
\be
\mu = \inf_{\psi \in E^+} \sup_{(x,t) \in D} \frac{L \psi}{\psi}
\ee
where $E^+ = \left\{ \psi \in E\;| \;\; \psi > 0 \right\}$.
\end{theorem}
{\bf Proof of Theorem \ref{th:ev1}:} This is an extension of the version proven in \cite{BH1} for time-independent operators. Clearly when $\psi = \phi$ is the eigenfunction for $L$, then $\frac{L \phi}{\phi} \equiv \mu$. Suppose there is a function $\psi > 0$, $\psi \in E$ such that
\be
\frac{L \psi}{\psi} < \mu  \;\;\;\text{for all} \;\;\;(x,t) \in D. \label{ev1}
\ee
Then $L \psi - \mu \psi = m < 0$ for some function m. Now there exists $\tau$ such that $\phi \leq \tau \psi$ for all $(x,t)$ with equality holding at a point. This follows from positivity and continuity of $\psi$ and $\phi$. Define the function $w \in E$ by $w = \frac{\phi}{\psi}$. Plugging $\phi = w \psi$ into the equation $L \phi - \mu \phi = 0$, we obtain
%\be
%\psi \left( a_{ij} w_{x_i x_j} - w_t + b_i w_{x_i}  \right ) + w ( L \psi - \mu \psi) + a_{ij} w_{x_i} %\psi_{x_j} + a_{ij} w_{x_j} \psi_{x_i} = 0. \no
%\ee
%Therefore,
\be
 a_{ij} w_{x_i x_j} - w_t + b_i w_{x_i}    + a_{ij} w_{x_i} \frac{\psi_{x_j}}{\psi} + a_{ij} w_{x_j} \frac{\psi_{x_i}}{\psi} = - \frac{w}{\psi} ( L \psi - \mu \psi) > 0 . \no
\ee
By definition of $w$ and $\tau$, $w \leq \tau $ with equality holding at one point. By maximum principle, we must have $w  \equiv \tau$. Note that periodicity in $t$ is necessary to apply the maximum principle in this way. Hence $\phi \equiv \tau \psi$, which contradicts (\ref{ev1}). Hence
\be
\sup_{(x,t) \in D} \frac{L \psi}{\psi} \geq \mu
\ee
for all $\psi \in E$. This completes the proof.

{\bf Proof of Lemma \ref{lem:ev1}:}
The Lemma extends an argument in \cite{BH1}, Section 5.2; we include the details for convenience. From the theorem we have
\be
\mu(\lambda) = \inf_{\psi \in E^+} \sup_{(x,t) \in D} \frac{L_\lambda \psi}{\psi}. \label{ev5}
\ee

For a function $\psi \in E^+$, we define $\hat \psi = e^{-\lambda k\cdot x} \psi$. Then it is easy to check that
\be
\frac{L_\lambda \psi }{ \psi} = \frac{\hat L \hat \psi}{\hat \psi} + f'(0) \no
\ee
where $\hat L \hat \psi = \Delta \hat \psi - \hat \psi_t + b \cdot \nabla \hat \psi$. Therefore, we can express $\mu(\lambda)$ as
\be
\mu(\lambda) =  f'(0) + \inf_{\psi \in \hat E^+_\lambda} \sup_{(x,t) \in D} \frac{\hat L \psi}{\psi}  \label{ev2}
\ee
with the set $\hat E^+_\lambda$ defined by $\hat E^+_\lambda = \left\{ e^{-\lambda k\cdot x} \psi \; | \;\; \psi \in E^+ \right\}$. For $\lambda, \gamma > 0$, let $\phi \in \hat E^+_\lambda$ and $\psi \in \hat E^+_\gamma$. If we define $w = \sqrt{\phi \psi}$, then $w \in E^+_{\alpha}$, where $\alpha = \frac{\lambda + \gamma}{2}$. Therefore, by (\ref{ev2})
\be
\mu(\frac{\lambda + \gamma}{2}) \leq f'(0) + \sup_{(x,t) \in D} \frac{\hat L w}{w}. \label{evw}
\ee
Let us compute
%\br
% \hat L  w &=& \frac{1}{2w} \left( \psi (\Delta \phi - \phi_t + b\cdot \nabla \phi)  \right) \no \\
% & & + \frac{1}{2w} \left( \phi (\Delta \psi - \psi_t + b\cdot \nabla \psi)  \right) \no \\
% & & - \frac{1}{4w^3} \left( \abs{\nabla \phi}^2 \psi^2 + \abs{\nabla \psi}^2 \phi^2 + 2 (\nabla \phi %\cdot \nabla \psi) \phi \psi  \right) \no \\
% & & + \frac{1}{w} \nabla \phi \cdot \nabla \psi. \no
%\er
%Therefore,
\br
 \frac{\hat L  w}{w} &=& \frac{1}{2}\frac{\hat L \phi}{\phi} + \frac{1}{2}\frac{\hat L \psi}{\psi}  \no \\
 & & - \frac{1}{4} \left( \frac{\abs{\nabla \phi}^2 }{\phi^2} + \frac{\abs{\nabla \psi}^2 }{\psi^2}  + 2\frac{ \nabla \phi \cdot \nabla \psi}{ \phi \psi}  \right) + \frac{\nabla \phi \cdot \nabla \psi}{\phi \psi} .\label{ev3}
\er
Note that
\br
2 \frac{\nabla \phi \cdot \nabla \psi}{\phi \psi} - \frac{\abs{\nabla \phi}^2 }{\phi^2} - \frac{\abs{\nabla \psi}^2 }{\psi^2}  = \no \\
 \frac{1}{\phi^2 \psi^2} \left( 2(\phi \psi) \nabla \phi \cdot \nabla \psi - \psi^2\abs{\nabla \phi}^2  - \psi^2 \abs{\nabla \psi}^2  \right). \label{ev4}
\er
If we take $a = \phi_{x_i} \psi$ and $b = \psi_{x_i} \phi$ and apply the inequality $2ab \leq a^2 + b^2$, we conclude that the expression in (\ref{ev4}) is nonpositive. Now returning to (\ref{ev3}), we conclude that
\be
\frac{\hat L  w}{w} \leq \frac{1}{2}\frac{\hat L \phi}{\phi} + \frac{1}{2}\frac{\hat L \psi}{\psi} . \no
\ee
Therefore,
\be
\sup_{(x,t) \in D} \frac{\hat L  w}{w} \leq  \frac{1}{2} \sup \frac{\hat L \phi}{\phi} + \frac{1}{2} \sup \frac{\hat L \psi}{\psi}.
\ee
Since $\phi$ and $\psi$ were arbitrarily chosen, we conclude from (\ref{ev2}) and (\ref{evw}) that
\be
\mu(\frac{\lambda + \gamma}{2}) \leq \frac{\mu(\lambda) + \mu(\gamma)}{2} .
\ee
Therefore, $\lambda \mapsto \mu(\lambda)$ is convex. This completes the proof.

{\bf Remark:} It follows easily from (\ref{ev2}) that if $\mu_\delta(\lambda)$ denotes principal eigenvalue of (\ref{ev6}) with $f'(0)$ replaced by $f'(0) - \delta$ for $\delta > f'(0) > 0$, we must have
\begin{itemize}
\item $\mu_\delta(\lambda) = \mu(\lambda) - \delta$ for all $\lambda > 0$
\item $\frac{\mu_\delta(\lambda)}{\lambda} \nearrow \frac{\mu(\lambda)}{\lambda}$ locally uniformly as $\delta \to 0$
\item $\frac{\mu_\delta(\lambda)}{\lambda}$ also satisfies the properties described in Lemma \ref{lem:ev1} and Proposition \ref{prop:ev1}
\end{itemize}

\section{Conclusions}
We proved the existence of KPP traveling fronts propagating at minimal speeds
in space-time periodic mean zero incompressble flows. The minimal speeds are
given by a variational principle defined through principal eigenvalue of
a linear space-time periodic parabolic operator. Front existence
holds also for general positive nonlinearity. It is interesting to
study of the qualitative
properties of minimal speeds as the advection field varies its intensity in
future works.

\section{Acknowledgments}
J.N. is grateful for support through an NSF VIGRE graduate fellowship at UT
Austin. M.R is supported by an NSF VIGRE Instructorship at UT Austin.
J.X is partially supported by NSF grant ITR-0219004.

\bibliographystyle{siam}

\end{document}